\documentclass{amsart}

\usepackage{amssymb}
\newcommand\Zs{{\mathbb Z}}
\newcommand\Cs{{\mathbb C}}

\begin{document}

\title{On the sigma function identity}

\author{Alexey~Gavrilov}
\email{gavrilov@lapasrv.sscc.ru}

\address{Institute of Computational Mathematics and Mathematical
Geophysics;~Russia,~Novosibirsk}

\begin{abstract}
We consider the known functional identity on the Weierstrass sigma function.
A complete classification of odd entire functions which satisfy the same
identity is obtained.
\end{abstract}

\maketitle

\section{Introduction}

Let $\Lambda$ be a lattice in $\Cs$ and
$\Lambda^{\prime}=\Lambda-\{0\}.$ Let
$$\sigma(z)=\sigma(z,\Lambda)=z\prod_{\lambda\in \Lambda^{\prime}}
\bigg{(}1-\frac{z}{\lambda}\bigg{)}e^
{\frac{z}{\lambda}+\frac{1}{2}(\frac{z}{\lambda})^2}\eqno{(1)}$$
be a Weierstrass sigma function. It is an odd entire
quasiperiodic function with $\Lambda$ as the set of zeros [1].

We shall deal with the identity
$$\sigma(x)\sigma(y)\sigma(z)\sigma(w)-\sigma(\frac{x+y+z-w}{2})
\sigma(\frac{x+y-z+w}{2})\sigma(\frac{x-y+z+w}{2})\sigma(\frac{-x+y+z+w}{2})-$$
$$-\sigma(\frac{x+y+z+w}{2})\sigma(\frac{x+y-z-w}{2})
\sigma(\frac{x-y+z-w}{2})\sigma(\frac{x-y-z+w}{2})=0, \eqno{(2)}$$
which holds for any $x,y,z,w\in \Cs.$
This identity may be easily derived from the classical Weierstrass parallelogramm
formula
$$\wp(x)-\wp(y)=-\frac{\sigma(x-y)\sigma(x+y)}{\sigma(x)^2\sigma(y)^2},$$
or from known properties of quasiperiodic functions or from the Riemann
theta identityes.
It is much more difficult to find it in the literature than to prove it.
The author knows only two papers where it has been written
(at least in an explicit form) [2,3]. In [2, lemma 1] it appears
as a new result. In [3] a reference to a book printed in 1893 was given.
Unfortunately the author has never read this book.

The aim of this paper is to prove the following

{\bf Theorem }
{\it Let  $f\neq 0$ be an odd entire function. Let
$$f(x)f(y)f(z)f(w)-f(\frac{x+y+z-w}{2})
f(\frac{x+y-z+w}{2})f(\frac{x-y+z+w}{2})f(\frac{-x+y+z+w}{2})-$$
$$-f(\frac{x+y+z+w}{2})f(\frac{x+y-z-w}{2})
f(\frac{x-y+z-w}{2})f(\frac{x-y-z+w}{2})=0 \eqno{(3)}$$
for any complex numbers $x,y,z,w.$ Then there exist $\alpha,\beta\in \Cs$
such that one of the following statements holds
\begin{eqnarray*}
\hbox{(1)} && \quad f(z)=ze^{\alpha  z^2+\beta};\\
\hbox{(2)} && \quad f(z)=\sin(az)e^{\alpha  z^2+\beta}\ \hbox{for some}\
   a\neq 0,\\
\hbox{(3)} && \quad f(z)=\sigma(z,\Lambda)e^{\alpha  z^2+\beta}\
\hbox{for some lattice}\ \Lambda.
\end{eqnarray*}
}

This result may be viewed as an algebraic definition of the
Weierstrass sigma function. It is clear that the statament of the theorem
is true under weaker hypothesis and an interesting question is which of the
conditions may be dropped.

\section{Some invariant}

The method of proof is the approximation of the function given by the function
of a required kind.
Let
${\mathcal M}=\{f: f(0)=0\}$ be the ideal in the algebra of entire functions
and let $\Omega$ be the set of odd entire functions which are not in
${\mathcal M}^2.$  For each $f\in\Omega$
we have
$$f(z)\equiv a_1 z+a_3 z^3+a_5 z^5+a_7 z^7 \mathop{mod}\, {\mathcal M}^9,$$
for some complex numbers $a_1, a_3, a_5, a_7,$ where $a_1\neq 0.$
Let
$$p(f)=a_3^2-2a_1a_5,$$
$$q(f)=3a_1^2a_7-3a_1a_3a_5+a_3^3.$$

{\bf Lemma 1 }
{\it  If $f\in\Omega$ and $p(f)=q(f)=0$ then there exist
$\alpha,\beta\in  \Cs$ such that
$$f(z)\equiv ze^{\alpha z^2+\beta} \mathop{mod}\, {\mathcal M}^9.$$}

The proof is a straightforward calculation.

If either of $p$ and $q$ is not equal to zero, then we may define
$$\mu(f)=\frac{p(f)^3}{q(f)^2}\in \Cs P^1.$$
We shall write $\mu(f(z))$ instead of $\mu(z\mapsto f(z)).$

{\bf Lemma 2 }
{\it Let $f_1,f_2\in\Omega$ and $\mu(f_1)=\mu(f_2)$ are well-defined.
Then
there exist $\alpha,\beta\in \Cs$ and $a\in  \Cs^{\times}$
such that
$$f_2(z)\equiv f_1(az)e^{\alpha z^2+\beta} \mathop{mod}\, {\mathcal M}^9.$$}

It is easy to see that
$$\mu(f(az)e^{\alpha z^2+\beta})=\mu(f(z))$$
for any $f\in\Omega.$
Let us choose $\alpha_i,\beta_i,\,i=1,2$ such that
$$\hat{f}_i(z)=f_i(z)e^{\alpha_i z^2+\beta_i}\equiv
z+A_iz^5+B_iz^7 \mathop{mod}\,  {\mathcal M}^9.$$

Then $p(\hat{f}_i)=-2A_i,\,q(\hat{f}_i)=3B_i.$
Since $\mu(\hat{f}_1)=\mu(\hat{f}_2),$ we have
$$\hat{f}_2(z)=\hat{f}_1(az)$$
for some $a\neq 0.$

\section{Modular forms}

We shall compute $\mu(\sigma(z,\Lambda)).$ For some reasons it is more
convenient to use the Jacobi function in this case instead of
the Weierstrass one. It is known that for the lattice
$$\Lambda=\rho(\Zs+\tau \Zs),\,\rho\neq 0,\, \Im(\tau)>0 $$
we have
$$\sigma(z,\Lambda)=\vartheta_1(\frac{z}{\rho}, \tau)e^{\alpha z^2+\beta},$$
where
$$\vartheta_1(z, \tau)=2\sum_{n=0}^{\infty}(-1)^ne^{\pi i\tau(n+\frac{1}{2})^2}
\sin((2n+1)\pi z)$$
is the first Jacobi theta function and $\alpha,\,\beta$ are some complex
numbers [1, \S 2.3]. Hence
$\mu(\sigma(z,\Lambda))=\mu(\vartheta_1(z, \tau)).$

Let
$$p(\tau)=p(\vartheta_1(z, \tau)),\, q(\tau)=q(\vartheta_1(z, \tau)).$$

{\bf Lemma 3 }
{\it The functions $p(\tau)$ and $q(\tau)$ may be written in the form
$$p(\tau)=\frac{\pi^2}{30}\eta^6g_2,\eqno{(4)}$$
$$q(\tau)=-\frac{\pi^3}{35}\eta^9g_3,\eqno{(5)}$$
where
$$\eta=\eta(\tau)=e^{\frac{\pi i\tau}{12}}
\prod_{n=1}^{\infty}(1-e^{2\pi in\tau})$$
is the Dedekind eta function and
$$g_2=g_2(\tau)=(2\pi)^4\bigg{(}\frac{1}{12}+20
\sum_{n=1}^{\infty}\frac{n^3}{e^{-2\pi in\tau}-1}\bigg{)},$$
$$g_3=g_3(\tau)=(2\pi)^6\bigg{(}\frac{1}{216}-\frac{7}{3}
\sum_{n=1}^{\infty}\frac{n^5}{e^{-2\pi in\tau}-1}\bigg{)}$$
are the Weierstrass modular forms.}

From the known equations
$$\vartheta_1(z, \tau+1)=e^{\frac{\pi i}{4}}\vartheta_1(z, \tau),$$
$$\vartheta_1(\frac{z}{\tau},-\frac{1}{\tau})=-i\sqrt{\frac{\tau}{i}}
e^{\frac{\pi iz^2}{\tau}}\vartheta_1(z, \tau)$$
we can derive
$$p(\tau+1)=ip(\tau),$$
$$p(-\frac{1}{\tau})=i\tau^7p(\tau).$$
The product $\eta^6g_2$ obeys the same equations, so
$\gamma(\tau)=\frac{p}{\eta^6g_2}$ is a modular function.
The direct computation shows that $\gamma(\tau)\to\frac{\pi^2}{30}$
as $\Im(\tau)\to\infty.$ The only zeros of the denominator are
$\tau$ which belong to the orbit of $\frac{-1+\sqrt{-3}}{2}$
under the modular group action. But this are first order zeros
whereas the order of any modular function at this points is divisible
by 3. Then $\gamma$ is bounded in upper halfplane so it is a constant.
The equality (5) may be proved by the same way.

As a corollary we have
$$\mu(\sigma(z,\Lambda))=\frac{p(\tau)^3}{q(\tau)^2}=
\frac{49}{1080}\frac{g_2^3}{g_3^2}=\frac{49}{40}\frac{j(\tau)}{j(\tau)-1728},$$
where
$$j(\tau)=e^{-2\pi i\tau}+744+196884 e^{2\pi i\tau}+\dots$$
is the modular invariant.

\section{Differential equation}

{\bf Lemma 4 }
{\it  Any odd entire function $f$ satisfying (3) obeys the equation
$$(f^{\prime}(0))^3f(2z)=f^4(z)(\ln f(z))^{\prime\prime\prime}. \eqno{(6)}$$}

Let us denote the left hand side of (3) by $F(x,y,z,w).$ Then
$$\frac{1}{6}\frac{d^3}{dt^3}F(x,x+t,x+\zeta t, x+\zeta^2 t)\bigg|_{t=0}=
-3f^2f^{\prime}f^{\prime\prime}+f^3f^{\prime\prime\prime}+2f(f^{\prime})^3-
(f^{\prime}(0))^3f(2z),$$
where$f=f(z)$ and $\zeta=e^{\frac{2\pi i}{3}}.$
This equals to zero, so (6) follows.

{\bf Lemma 5 }
{\it  Let $f_1,f_2\in\Omega$ be two functions satisfying (6).
If $f_2\equiv f_1 \mathop{mod}\, {\mathcal M}^9$ then
$f_2=f_1$ identically.}

We may asuume $f_1^{\prime}(0)=f_2^{\prime}(0)=1.$
If $f_2\neq f_1$ then for some odd $n\ge 9$
$$f_2(z)\equiv f_1(z)+bz^n \mathop{mod}\, {\mathcal M}^{n+2},$$
where $b\neq 0.$ We have
$$f_2(2z)-f_1(2z)=b2^nz^n \mathop{mod}\, {\mathcal M}^{n+2},$$
$$f_2^4(z)(\ln f_2(z))^{\prime\prime\prime}-
f_1^4(z)(\ln f_1(z))^{\prime\prime\prime}=
bz^n[(n-1)(n-2)(n-3)+8] \mathop{mod}\, {\mathcal M}^{n+1}.$$
Since $f_1$ and $f_2$ are both solutions of (6),
$$bz^n\psi(n)\equiv 0 \mathop{mod}\, {\mathcal M}^{n+1},$$
where $\psi(n)=(n-1)(n-2)(n-3)+8-2^n.$
By assumption $b\neq 0$ so $\psi(n)=0$ but this is imposible for $n\ge 9.$

Now we can prove the theorem. Let $f$ be a function which satisfies the
conditions. By Lemma 4 it is a solution of (6). So
$f\in\Omega$ because nonzero odd function may not be a solution of
$f^4(z)(\ln f(z))^{\prime\prime\prime}=0.$
Now it is enough to prove that there exists a function of required kind
approximating $f$ up to 9th order, then the conclusion of theorem
holds by Lemma 5.

If $p(f)=q(f)=0$ then this is the case by Lemma 1. If
$\mu(f)=\frac{49}{40}$ then $\mu(f)=\mu(\sin(z)).$ If
$\mu(f)\neq\frac{49}{40}$ then
$\mu(f)=\frac{49}{40}\frac{j(\tau)}{j(\tau)-1728}$
for some $\tau$ from the upper halfplane so
$$\mu(f)=\mu(\sigma(z,\Lambda)),\,\Lambda=\Zs+\tau \Zs$$
(the case $\mu=\infty$ is included).
In both cases $f$ has the required approximation by Lemma 2.
The proof is complete.

\section*{}
1. Hurwitz A., Courant R. Vorlesungen uber allgemeine
funktionentheorie und elliptische funktionen, II. Springer-Verlag, 1964.\\
2. McCullough S., Shen Li-Chien. On the Szego kernel of an annulus//
Proc. Amer. Math. Soc. 121(4), 1994, p. 1111-1122.\\
3. Amdeberhan T. A determinant of the Chudnovskys generalazing the
elliptic Frobenius-Stickelberger-Cauchy determinant identity//
Electron J. Comb. 7(1), 2000, N6.
\end{document}